\documentclass[11pt,letterpaper]{article}

\usepackage{fullpage}

\usepackage{amsmath,amsthm,amssymb}
\usepackage{mathtools}
\usepackage{mathrsfs}
\usepackage[shortlabels]{enumitem}

\usepackage{comment}

\usepackage{hyperref}
\usepackage[nameinlink,capitalise,noabbrev]{cleveref}

\usepackage{graphicx}
\usepackage{float}
\usepackage{tikz}

\theoremstyle{definition}
\newtheorem{definition}{Definition}

\newtheorem{remark}[definition]{Remark}

\theoremstyle{plain}

\newtheorem{theorem}[definition]{Theorem}
\newtheorem{lemma}[definition]{Lemma}
\newtheorem{proposition}[definition]{Proposition}
\newtheorem{claim}[definition]{Claim}

\newtheorem{corollary}[definition]{Corollary}

\def \l {\ell}

\def \ex {\mathrm{ex}}
\def \sm {\setminus}
\def \ce {\coloneqq}
\def \E {\mathbb{E}}

\renewcommand{\le}{\leqslant}
\renewcommand{\ge}{\geqslant}
\renewcommand{\leq}{\leqslant}

\def \eps {\varepsilon}
\def \es {\varnothing}
\renewcommand \b[2] {\binom{#1}{#2}}

\def \A {\mathcal{A}}

\def \C {\mathcal{C}}
\def \F {\mathcal{F}}
\def \G {\mathcal{G}}
\def \mH{\mathcal{H}}
\def \mP{\mathcal{P}}

\title{On the maximum $F$-free induced subgraphs \\ in $K_t$-free graphs}

\author
{
J{\'o}zsef Balogh\thanks{University of Illinois Urbana-Champaign, 1409 W. Green Street, Urbana IL 61801, United States.
E-mail: \texttt{\{jobal, cechen4, haoranl8\}@illinois.edu}.
}
\thanks{Research is supported in part by NSF grants DMS-1764123 and RTG DMS-1937241, FRG DMS-2152488, the Arnold O. Beckman Research Award (UIUC Campus Research Board RB 24012).}
\and
Ce Chen\footnotemark[1]
\and
Haoran Luo\footnotemark[1] \thanks{Research is partially supported by Trjitzinski Fellowship.}
}

\date{}

\begin{document}
\maketitle

\begin{abstract}
    For graphs $F$ and $H$, let $f_{F,H}(n)$ be the minimum possible size of a maximum $F$-free induced subgraph in an $n$-vertex $H$-free graph. This notion generalizes the Ramsey function and the Erd\H{o}s--Rogers function.
    Establishing a container lemma for the $F$-free subgraphs, we give a general upper bound on $f_{F,H}(n)$, assuming the existence of certain locally dense $H$-free graphs. In particular, we prove that for every graph $F$ with $\ex(m,F) = O(m^{1+\alpha})$, where $\alpha \in [0,1/2)$, we have
    \[
     f_{F, K_3}(n) =  O\left(n^{\frac{1}{2-\alpha}}\left(\log n\right)^{\frac{3}{2-
 \alpha}}\right)
    \quad
    \textrm{and}
    \quad
    f_{F, K_4}(n) =  O\left(n^{\frac{1}{3-2\alpha}}\left(\log n\right)^{\frac{6}{3-2\alpha}}\right).
    \]
For the cases where $F$ is a complete multipartite graph,
letting $s = \sum_{i=1}^r s_i$,
we prove that
    \[
        f_{K_{s_1,\ldots,s_r}, K_{r+2}}(n) = O \left( n^{\frac{2s -3}{4s -5}} (\log n)^{3} \right).
    \]
    We also make an observation which improves the bounds of $\ex(G(n,p),C_4)$ by a polylogarithmic factor.
\end{abstract}

\section{Introduction}

For graphs $F$ and $G$, we say that a subset $S \subseteq V(G)$ is \emph{$F$-independent} if the induced subgraph $G[S]$ does not contain $F$ as a subgraph. The \emph{$F$-independence number} $\alpha_F(G)$ of $G$ is
the maximum size of an $F$-independent subset in $G$.
For graphs $F$ and $H$, let $f_{F,H}(n)$ be the minimum possible value of $\alpha_F(G)$ over all $n$-vertex $H$-free graphs $G$. If $F= K_i$ and $H = K_j$, i.e., they are complete graphs, then we simply write $f_{i,j}(n)$ for $f_{K_i,K_j}(n)$.

In 1962, Erd\H{o}s and Rogers~\cite{erdos1962construction} initiated the study of $f_{s,t}(n)$, which is now commonly referred to as the Erd\H{o}s--Rogers function.
Note that the case $s=2$ is related to  Ramsey theory~\cite{ramsey1929problem}: the Ramsey number $r(s,t)$ is at least $n$ if and only if $f_{2, t}(n) < s$.
For example, by the bound $r(3,t) = \Theta(t^2/\log t)$~\cite{ajtai1980note, kim1995ramsey} and the recent breakthrough  $\Omega(t^3 / (\log t)^4) = r(4,t) = O(t^3 / (\log t)^2)$~\cite{ajtai1980note, mattheus2024asymptotics}, we have
\begin{equation} \label{equ::f23f24}
    f_{2,3}(n) = \Theta\left(n^{\frac{1}{2}} (\log n)^{\frac{1}{2}}\right)
\quad \textrm{and} \quad
   \Omega \left(n^{\frac{1}{3}} (\log n)^{\frac{2}{3}}\right) =  f_{2,4}(n) = O\left(n^{\frac{1}{3}} (\log n)^{\frac{4}{3}}\right).
\end{equation}

For $s > 2$, Erd\H{o}s and Rogers~\cite{erdos1962construction} proved that $f_{s,s+1}(n) \le n^{1-\Omega_s(1)}$.
In the following decades, there were many papers on this topic and significant improvements have been achieved for various pairs of $(s,t)$. The current best bound~\cite{dudek2014generalized,mubayi2024order} on $f_{s,s+1}(n)$ is
$$
    \Omega\left( (n \log n / \log\log n)^{\frac{1}{2}} \right) = f_{s,s+1}(n) = O\left(n^{\frac{1}{2}} \log n\right),
$$
and for $f_{s,s+2}(n)$, the best known bound~\cite{sudakov2005large, janzer2023improved} is
$$
    n^{\frac{1}{2} - \frac{1}{6s-6}} (\log n)^{\Omega(1)} = f_{s,s+2}(n) = O\left(n^{\frac{2s-3}{4s-5}}(\log n)^3\right).
$$
There are also many results on general $(s,t)$, see~\cite{bollobas1991graphs, krivelevich1994Ksfree, krivelevich1995bounding, sudakov2005new, sudakov2005large}.
For more results about the Erd\H{o}s--Rogers function, we refer the reader to the recent paper~\cite{janzer2023improved} and the citations within.

In this paper, we investigate some general cases of $f_{F,H}(n)$, where $F$ is not necessarily a complete graph. This problem was suggested by Dhruv Mubayi, when he visited UIUC in  Fall 2023.  By definition, we have $f_{F_1,H}(n) \le f_{F_2, H}(n)$ whenever $F_1 \subseteq F_2$, so $f_{2,H}(n)$ serves as a lower bound of $f_{F,H}(n)$ for every $F$ containing at least one edge. For the upper bound\footnote{We are grateful to Jacques Verstra{\"e}te for pointing out this to us.}, note that
every $F$-free subgraph in any graph $G$ has size
less than $r(F, \alpha(G) + 1)$. Hence,
every upper bound on $r(F,t)$ and $f_{2,H}(n)$ offers an upper bound for $f_{F,H}(n)$. For example, using~\eqref{equ::f23f24} and $r(C_{\l} , t) = O(t^{1+ 1 / (\lceil \l/2 \rceil -1)})$~\cite{erdos1978cycle}, we deduce that for every $k \ge 2$,
$$
\Omega\left(n^{\frac{1}{2}} (\log n)^{\frac{1}{2}}\right)    = f_{C_{2k}, K_3} (n) = O\left(n^{\frac{k}{2k-2}}\left(\log n\right)^{\frac{k}{2k-2}}\right),
$$ and
$$
\Omega \left(n^{\frac{1}{3}} (\log n)^{\frac{2}{3}}\right) = f_{C_{2k}, K_4} (n) =
O\left(n^{\frac{k}{3k-3}}\left(\log n\right)^{\frac{4k}{3k-3}}\right),
$$
and the same bounds hold for $f_{C_{2k-1}, K_3}(n)$ and $f_{C_{2k-1}, K_4}(n)$.

Our first result is an improvement on the upper bounds of $f_{F, K_3}(n)$ and $f_{F, K_4}(n)$ where $F$ has a small Tur\'an number.
Recall that
the \emph{Tur\'an number} $\ex(n,F)$ of $F$ is the maximum number of edges in an $n$-vertex $F$-free graph.
\begin{theorem} \label{thm::FK3K4}
For every graph $F$ with $\ex(b,F) \le O(b^{1+\alpha})$, where $\alpha \in [0,1/2)$, we have
    \[
     f_{F, K_3}(n) =  O\left(n^{\frac{1}{2-\alpha}}\left(\log n\right)^{\frac{3}{2-
 \alpha}}\right)
    \quad \quad
    \textrm{and}
    \quad \quad
    f_{F, K_4}(n) =  O\left(n^{\frac{1}{3-2\alpha}}\left(\log n\right)^{\frac{6}{3-2\alpha}}\right).
    \]
\end{theorem}
\noindent
In particular, by the Bondy--Simonovits theorem~\cite{bondy1974cycles}, we have $\ex(b, C_{2k}) = O(b^{1+1/k})$ for every integer $k \ge 2$. Therefore, we have the following corollary.
\begin{corollary}\label{thm::c2k}
	For every integer $k \ge 3$, we have
    \[
	f_{C_{2k}, K_3}(n) =  O\left(n^{\frac{k}{2k-1}}\left(\log n\right)^{\frac{3k}{2k-1}}\right)
    \quad \quad
    \textrm{and}
    \quad \quad
    f_{C_{2k}, K_4}(n) =  O\left(n^{\frac{k}{3k-2}}\left(\log n\right)^{\frac{6k}{3k-2}}\right).
    \]
\end{corollary}
\noindent
The main tool for the proof of \cref{thm::FK3K4} is a container lemma (\cref{lem::container}), which offers an upper bound on the number of $F$-independent sets of a certain size in a ``locally dense'' graph.
Using this Container Lemma, we can actually prove the following much more general result.
\begin{theorem} \label{thm::generalfFG}
    For real numbers $\alpha$, $\beta$, $\gamma$, $\delta$, and $\theta$, where $0 \le \alpha  < 1$, $0 <\beta, \delta, \gamma$, $0 <\theta < 1$, and $\beta < \theta(1-\alpha)$, and a graph $F$ with $\ex(b,F) \le O(b^{1+\alpha})$,
    there exist positive real numbers $c$, $C$, and $N$ such that the following holds for every integer $n > N$. Let
    $$
        s \ce
        \left \lceil n^{\frac{\beta}{1-\alpha}} (\log n)^{\frac{3}{1-\alpha}}
        \right \rceil
        .
    $$
    If $G$ is an $n$-vertex graph with $e(G[A]) \ge  |A|^2 \cdot \delta n^{-\beta}$ for every $A \subseteq V(G)$ with $|A| > \gamma n^{\theta}$, then there exists $G'\subseteq G$ with at least
    $
        m \ce c {s} n^{1-\theta}
    $
    vertices and
    $$
        \alpha_F(G') \le s \le C \cdot
        m^{\frac{\beta}{(1-\alpha)(1-\theta)+\beta}}
        (\log m)^{\frac{3(1-\theta)}{(1-\alpha)(1-\theta)+\beta}}.
    $$
\end{theorem}
\noindent
We will derive \cref{thm::FK3K4} from \cref{thm::generalfFG} by letting $G$ be the pseudorandom $K_3$-free graph given by Alon~\cite{alon1994explicit} and the construction of certain $K_4$-free graphs given by Mattheus and Verstra{\"e}te~\cite{mattheus2024asymptotics}, respectively.

\cref{thm::generalfFG} also implies many other results, such as the following theorem, which gives an upper bound on $f_{C_{2k},K_t}(n)$,
assuming the existence of ``perfectly'' pseudorandom $K_t$-free graphs.
Recall that an \emph{$(n,d, \lambda)$-graph} is an $n$-vertex $d$-regular graph such that the absolute value of every eigenvalue of its adjacency matrix, but the largest one, is at most $\lambda$.
For more background about $(n,d,\lambda)$-graphs, we refer the reader to the survey~\cite{Krivelevich2006pseudorandom}.

\begin{theorem} \label{thm::c2kpseu}
    For every graph $F$ and integer $t \ge 3$, where $\ex(b,F) = O(b^{1+\alpha})$ for some $\alpha \in [0, (t-2)/(t-1))$,
    if for every sufficiently large $m$, there exists a $K_t$-free $(m, d, \lambda)$-graph $G_m$ with $d = \Theta(m^{1- 1/(2t-3)})$ and $\lambda=O(d^{1/2})$, then we have
    $$
        f_{F, K_t}(n) = O\left( n^{\frac{1}{(1-\alpha)(t-2)+1}} (\log n)^{\frac{3(t-2)}{(1-\alpha)(t-2)+1}} \right).
    $$
\end{theorem}
\noindent
An interesting observation is that, while different $K_4$-free graphs are used, Theorems~\ref{thm::FK3K4} and~\ref{thm::c2kpseu} give the same upper bound on $f_{F,K_4}(n)$ for every graph $F$ with $\ex(b,F) = O(b^{1+\alpha})$ for some $\alpha \in [0,1/2)$.

For $F = C_4$ and $H = K_4$, Theorems~\ref{thm::generalfFG} and~\ref{thm::c2kpseu} fail to provide a reasonable bound on $f_{F,H}(n)$.
Instead, we have the following slightly more general theorem. Denote by $K_{s_1, \ldots, s_r}$ the complete $r$-partite graph with parts of sizes $s_1, \ldots, s_r$, respectively.
\begin{theorem} \label{thm::Ks1srVsKr+2}
    For every integer $r \ge 2$ and positive integers $s_1, \ldots, s_r$, letting $s \ce \sum_{i=1}^r s_i$,
    we have
    $$
    f_{K_{s_1,\ldots,s_r}, K_{r+2}}(n) = O \left( n^{\frac{2s -3}{4s -5}} (\log n)^{3} \right).
    $$
\end{theorem}
\noindent
In particular, \cref{thm::Ks1srVsKr+2} implies
$$
f_{C_4, K_4} (n) = O\left(n^{\frac{5}{11}} (\log n)^{3}\right)
\quad
\text{and in general }
\quad
f_{K_{a,a}, K_4} (n) = O\left(n^{\frac{1}{2}-\frac{1}{16a - 10}} (\log n)^{3}\right).
$$
Our proof for \cref{thm::Ks1srVsKr+2} follows closely the recent proof of $f_{{r}, {r+2}}(n) = O(n^{(2r-3)/(4r-5)}(\log n)^3)$ by Janzer and Sudakov~\cite{janzer2023improved}, which followed closely the work of Mattheus and Verstra{\"e}te~\cite{mattheus2024asymptotics}. We generalize their ideas to bound the number of copies of $K_{s_1,\ldots, s_r}$ in the graphs based on the construction of $K_4$-free graphs given in~\cite{mattheus2024asymptotics}.

\begin{remark}
    While working on this project, we learned that Mubayi and Verstra{\"e}te~\cite{mubayi2024Erdos} proved that $f_{F,K_3}(n) = n^{1/2+o(1)}$ for every triangle-free graph $F$, and for every $K_4$-free graph $F$ containing a cycle, there is a constant $c_F > 0$ such that $f_{F,K_4}(n) > n^{1/3+ c_F + o(1)}$.
    In particular, their results imply (and somewhat show the tightness of some of our results) that $f_{C_{2k},K_4} (n) =n^{1/3+\Theta(1/k)}$.
\end{remark}

We also have the following observation on the random Tur\'an number of $C_{4}$.  Recall that $\ex(G(n,p),C_4)$ is the maximum number of edges in a $C_4$-free subgraph of the Erd\H{o}s--R\'enyi graph $G(n,p)$.

\begin{proposition} \label{pro::impBoundGnpC4}
If $p \ge n^{-1/3} (\log n)^{8/3}$, then with high probability, we have
$$
 \ex(G(n,p), C_4)= O \left(p^{1/2}n^{3/2}\right).
$$
\end{proposition}
\noindent
\cref{pro::impBoundGnpC4} improves a result of Morris and Saxton~\cite{morris2016number}, which gives the same upper bound under the assumption that $p \ge n^{-1/3} (\log n)^{4}$.
See \cref{sec::GnpC4} for more discussion on this problem.

The rest of this paper is organized as follows. In \cref{sec::container}, we prove our Container Lemma and use it to prove Theorems~\ref{thm::generalfFG},~\ref{thm::FK3K4}, and~\ref{thm::c2kpseu}. In \cref{sec::boundByJanzerSudakov}, we prove \cref{thm::Ks1srVsKr+2}, and in \cref{sec::GnpC4}, we prove \cref{pro::impBoundGnpC4}.

\textbf{Notation.} For a graph $G$, a subset $S \subseteq V(G)$, and a vertex $v \in V(G)$, we write $N_S(v)$ for the set of \emph{neighbors} of $v$ in $S$ and let $d_S(v) \ce |N_S(v)|$ be the \emph{degree} of $v$ in $S$. When $S = V(G)$, we write $N_G(v)$ for $N_{V(G)}(v)$ and $d_G(v)$ for $d_{V(G)}(v)$. All logarithms are in base $e$, except in \cref{sec::boundByJanzerSudakov}, where they are in base $2$.

\section{The Container Method}
Our main tool for proving \cref{thm::generalfFG} is the following Container Lemma, whose proof is an extended variant of the graph container method. See~\cite{balogh2020number} for a prior similar application of the method.
\label{sec::container}
\begin{lemma}[Container Lemma] \label{lem::container}
For every graph $F$ and positive real numbers $d, r$, if
$G$ is an $n$-vertex graph, where $e(G[A]) \ge d \cdot |A|^2$
for every $A \subseteq V(G)$ with $|A| > r$,
then for every $F$-independent subset $S \subseteq V(G)$ with $\ex(|S|, F) \cdot d \ge \log n$, there exist $S_h$, $T$, $C \subseteq V(G)$, where $C$ depends only on $S_h$, $T$, and $|S|$, such that
\begin{enumerate}[(i)]
     \item $(S_h \cup T) \subseteq S$ and $S \sm (S_h \cup T) \subseteq C$,
     \item $|S_h| \le 2 \cdot \left( \frac{\ex(|S|,\,F) \log n}{d} \right)^{1/2}$,
     \item $|T|\le 32 \cdot \left( \frac{\ex(|S|,\,F) \log n}{d} \right)^{1/2}$,
     \item $|C|\le r$.
\end{enumerate}
\end{lemma}

\begin{proof}
    Choose $S \subseteq V(G)$ such that $G[S]$ is $F$-free and $\ex(|S|, F) \cdot d \ge \log n$. Letting $s \ce |S|$, we have $e(G[S]) \le \ex(s, F)$. We use the following algorithm which takes $S$ and outputs the sets $S_h$, $T$, and $C$.

    Let $\Delta \ce (\ex(s,F) \cdot d / \log n)^{1/2} \ge 1$.
    Let
    $
    S_h \ce \{
    v \in S: d_S(v) \ge \Delta \},
    $
    and $S_\ell \ce S \sm S_h$. Fix an arbitrary ordering $\pi$ of $V(G)$.
    We start the algorithm with sets $T_0 \ce \es$, $A_0 \ce V(G) \sm S_h$, and the function $t_0(v) \ce 0$ for every vertex $v\in V(G)$.
    For every integer $i \ge 0$, we pick a vertex $u_i \in A_i$ of the maximum degree in $G[A_i]$. In case there is more than one such vertex, we choose the one that comes earlier in $\pi$.

    If $u_i\in S$, then let $T_{i+1}\ce T_i\cup\{u_i\}$,
    $$
    t_{i+1}(v) \ce \left\{
    \begin{array}{ll}
        t_i(v)+1 & \textrm{if\ } v \in N_{G[A_i]}(u_i),\\
		t_i(v)   & \textrm{otherwise,}
    \end{array}
    \right.
    $$
    and
    $$
        A_{i+1} \ce A_i \sm (\{u_i\} \cup  \{v\in A_i: t_{i+1}(v)>  \Delta\}).
    $$
    If $u_i \notin S$, then let $T_{i+1}\ce T_i$, $t_{i+1}(v)\ce t_i(v)$ for every vertex $v\in V(G)$, and $A_{i+1}\ce A_i\sm\{u_i\}$. The algorithm terminates at step $K$, where $K$ equals the smallest $i$ such that $|A_i| \le r$, and outputs
    $$
    S_h, \quad
	T\ce T_K, \quad \text{and} \quad C \ce A_K.
    $$

    Note that $|A_{i+1}| \le |A_i| -1$ for every $i \ge 0$. Hence, the algorithm above indeed terminates after finitely many steps.
    Now, we verify that the outputs $S_h$, $T$, and $C$ satisfy all the claimed properties.

    \begin{claim}
        For inputs $S'$, $S'' \subseteq V(G)$ where $|S'| = |S''|$, if the algorithm outputs the same $S_h$ and $T$, then the algorithm outputs the same $C$.
    \end{claim}
    \begin{proof}
        We write $\Delta', S'_h, S'_\ell, T'_i, A'_i, t'_i, u'_i, K', T', C'$ and $\Delta'', S''_h, S''_\ell, T''_i, A''_i, t''_i, u''_i, K'', T'', C''$ for the corresponding parameters in the algorithm with input $S'$ and $S''$, respectively. By assumption, we have $S'_h = S''_h$ and hence, $A'_0 = A''_0$ and $u'_0 = u''_0$. Since $|S'| = |S''|$, we also have $\Delta' = \Delta''$.
        Besides, we have $T'_0 = T''_0 =\es $ and $t'_0 = t''_0$ by definition.
        Now, observe that if $C' \neq C''$, then there is an index $i$ such that $u'_i = u''_i$ but $u'_i \in S'$, $u''_i \notin S''$. Then, we have $u'_i \in T'$ but $u''_i \notin T''$, a contradiction to the assumption that $T' =T''$.
    \end{proof}

    \begin{claim}
        $(S_h \cup T) \subseteq S$ and $S \sm (S_h \cup T) \subseteq C$.
    \end{claim}
    \begin{proof}
        By definition, we have $S_h \subseteq S$ and $T = T_K \subseteq S_\ell \subseteq S$, so $(S_h \cup T) \subseteq S$.
        For every vertex $v \in S \sm (S_h \cup T) = S_\ell \sm T_K$, if $v \notin C =A_K$, then there is an index $i < K$ such that $t_{i+1}(v) > \Delta$. This means that $v$ has more than $\Delta$ neighbors in $\{u_1,\ldots, u_i\} \cap S$, so $v \in S_h$, a contradiction.
    \end{proof}

    \begin{claim}
        $|S_h| \le 2 \cdot \left( \frac{\ex(s,\,F) \log n}{d} \right)^{1/2}$.
    \end{claim}
    \begin{proof}
        We have
        $$
            |S_h| \cdot \Delta \le \sum_{v \in S} d_S(v) = 2 e(G[S]) \le 2 \ex(s, F).
        $$
        Recalling that $\Delta = (\ex(s,F) \cdot d / \log n)^{1/2}$, we have
        \begin{equation*}
          |S_h| \le \frac{2\cdot \ex(s,F)}{\Delta} = 2 \left( \frac{\ex(s,\,F) \log n}{d} \right)^{1/2}. \qedhere
        \end{equation*}
    \end{proof}

    \begin{claim}
        $|T|\le 32 \cdot \left( \frac{\ex(s,\,F) \log n}{d} \right)^{1/2}$.
    \end{claim}
    \begin{proof}
    For every positive integer $j \in [1, 2\log n]$, let $I_j \ce \{i \in [0,K-1]: n/2^j < |A_{i}| \le n/2^{j-1}\}$.
    Let $J \ce \{j: I_j \neq \es\}$.
    Note that  $\{I_j: j \in J\}$ forms a partition of $[0, K-1]$, as $2^{2\log n +1 -1} > n$. For every $j\in J$, let $a_j$ and $b_j$ be the minimum and maximum elements in $I_j$, respectively.
    Since $T_0 = \es$, we have that $T_K$ can be expressed as the disjoint union
    \[
        T_K =
        \bigsqcup_{j \in J}
        \left(
        \left(T_{b_j} \sm T_{a_j}\right)
        \cup \left( T_{b_j+1} \sm T_{b_j}\right)
        \right),
    \]
    so
    \begin{equation} \label{equ::disUniTK}
        |T| = |T_K| \le 2 \log n + \sum_{j\in J} |T_{b_{j}} \sm T_{a_{j}}|.
    \end{equation}
    For every $j \in J$, we bound $|T_{b_{j}} \sm T_{a_{j}}|$ by the observation
    \begin{equation}\label{equ::conti_d_t}
        \sum_{\substack{i\,:\,a_j \le i < b_j \\ u_i \in T_K}} d_{G[A_i]}(u_i)
        \le
        \sum_{v \in A_{a_j}} t_{b_j} (v),
    \end{equation}
    which holds because $t_{i+1}(v)$ increases by one for every $v \in N_{G[A_i]}(u_i) \subseteq A_{a_i}$ when moving $u_i$ into $T_i$. For the left-hand side of~\eqref{equ::conti_d_t}, recall that $u_i$ is the vertex with maximum degree in $G[A_i]$, so
    $$
        d_{G[A_i]}(u_i) \ge \frac{2 e(G[A_i])}{ |A_i|} \ge \frac{2 e(G[A_i])}{|A_{a_j}|}.
    $$
    Note that for every $i < K$, we have $|A_i| > r$ by the terminating condition of our algorithm. By the assumption on $G$, we have $e(G[A_i]) \ge d\cdot |A_i|^2 \ge d \cdot |A_{b_j}|^2$. Hence,
    $$
        d_{G[A_i]}(u_i) \ge \frac{2d|A_{b_j}|^2}{|A_{a_j}|}
        \ge \frac{2|A_{b_j}|}{|A_{a_j}|} \cdot d|A_{b_j}|
        \ge d |A_{b_j}|
        .
    $$
    For the right-hand side of~\eqref{equ::conti_d_t}, note that $t_i(v)$ is never larger than $\Delta + 1$ for any $i$ and $v \in V(G)$. Therefore, by~\eqref{equ::conti_d_t},
    $$
        |T_{b_j} \sm T_{a_j}| \cdot d|A_{b_j}| \le |A_{a_j}| \cdot (\Delta + 1),
    $$
    which implies that
    $$
        |T_{b_j} \sm T_{a_j}| \le \frac{(\Delta + 1)}{ d} \cdot \frac{|A_{a_j}|}{|A_{b_j}|} \le \frac{2 (\Delta+1) }{d}
        \le \frac{4 \Delta}{d}
        .
    $$
    Note that $\log n \le (\ex(s,\,F) \log n/d)^{1/2}$ by our assumption $\ex(s, F) \cdot d \ge \log n$.
    By~\eqref{equ::disUniTK}, we have
    \begin{equation*}
        |T| = |T_K| \le
        2 \log n + \frac{4\Delta}{d} \cdot 2 \log n
        \le 32 \left( \frac{\ex(s,\,F) \log n}{d} \right)^{1/2}. \qedhere
    \end{equation*}
    \end{proof}

    \begin{claim}
        $|C| \le r$.
    \end{claim}
    \begin{proof}
        The algorithm terminates only when $|C| = |A_K| \le r$.
    \end{proof}

    Therefore, the sets $S_h$, $T$, and $C$ indeed satisfy all the claimed properties.
\end{proof}

    Now, we prove \cref{thm::generalfFG} and use it to prove Theorems~\ref{thm::FK3K4} and~\ref{thm::c2kpseu}.
    \begin{proof}[Proof of \cref{thm::generalfFG}]
        Let $c > 0$ be sufficiently small and $C, N$ be sufficiently large. Without loss of generality, we may assume that $\ex(b, F) = \Omega(b)$. Recalling that
        $
            s =
        \left \lceil n^{\frac{\beta}{1-\alpha}} (\log n)^{\frac{3}{1-\alpha}} \right \rceil
        $
        and $0\le \alpha < 1$,
        we have
        $$
        \ex\left(s, F\right) \cdot \delta n^{-\beta}
        = \Omega\left(s \cdot  n^{-\beta}\right)
        = \Omega\left( n^{\frac{\beta}{1-\alpha}-\beta} (\log n)^{\frac{3}{1-\alpha}} \right) > \log n.
        $$
        Hence, we can apply \cref{lem::container} to the $F$-independent subsets $S \subseteq V(G)$ of size $s$
        with parameters $d = \delta n^{-\beta}$ and $r = \gamma n^{\theta}$.
        Let
        $$
            t \ce \left\lfloor 40 \cdot \left({\ex(s ,\,F) \log n} / {(\delta n^{-\beta})} \right)^{1/2} \right\rfloor,
        $$
        which is an upper bound for both
        $|S_h|$ and $|T|$ in \cref{lem::container}, then
        \begin{equation} \label{equ::tlogn}
            t \log n
            = O \left( \left(n^{\frac{\beta}{1-\alpha}\cdot (1+\alpha) + \beta} (\log n)^{\frac{3}{1-\alpha} (1+\alpha)+1} \right)^{\frac{1}{2}}\right)\cdot \log n
            = O\left( n^{\frac{\beta}{1-\alpha}} (\log n)^{\frac{2+\alpha}{1-\alpha}} \right) \cdot \log n = O(s)
            .
        \end{equation}
        Note that $s= o(n^{\theta})$ by the assumption that $\beta < \theta(1-\alpha)$.
        By \cref{lem::container}, the number of $F$-independent subsets of $V(G)$ of size $s$ is at most
        $$
            \sum_{i=0}^{t}\sum_{j=0}^{t} \b{n}{i} \b{n}{j} \b{\gamma n^{\theta}}{s}
            \le (t+1)^2 n^{2t} \b{\gamma n^{\theta}}{s}
            \le n^{4t} \left(\frac{e\gamma n^{\theta}}{s}\right)^{s}
            \le e^{4t\log n} \left(\frac{e\gamma n^{\theta}}{s}\right)^{s}
            .
        $$

        Let $p \ce 2c\cdot s/n^{\theta}$. We choose a random subset $B$ of $V(G)$, where every vertex is included with probability $p$, independently of each other. Let $\F$ be the (random) collection of $F$-independent subsets of $B$ of size $s$.
        Then, we have $\E[|B|] = np$ and
        $$
            \E[|\F|] \le e^{4t\log n}\left(\frac{e\gamma n^{\theta}}{s}\right)^{s} p^{s}
            \le e^{4t\log n} \left(2ce\gamma \right)^{s}
            \stackrel{\eqref{equ::tlogn}}{<} 1.
        $$
        Let $B'$ be the set obtained from $B$ by removing one vertex from every $S \in \F$ and $G' \ce G[B']$. Then, there is a choice of $G'$ such that $G'$ has no $F$-free subgraph of size $s$ and
        \begin{equation*}
            |V(G')| \ge \E[|B'|] \ge \E[|B| - |\F|] = \E[|B|] - \E[|\F|]  \ge np -1
            =2csn^{1-\theta} -1
            \ge csn^{1-\theta}.
        \end{equation*}
        Note that for $m \ce csn^{1-\theta}$, we have $(1-\theta) \log n \le \log m \le \log n$ and
        \[
        m = csn^{1-\theta} \ge cn^{\frac{\beta}{1-\alpha}}(\log n)^{\frac{3}{1-\alpha}} \cdot n^{1-\theta}
        = cn^{\frac{(1-\alpha)(1-\theta) + \beta}{1-\alpha}} (\log n)^{\frac{3}{1-\alpha}},
        \]
        so
        \[
            s \le 2n^{\frac{\beta}{1-\alpha}} (\log n)^{\frac{3}{1-\alpha}} \le C\cdot m^{\frac{\beta}{(1-\alpha)(1-\theta)+\beta}} (\log m)^{\frac{3(1-\theta)}{(1-\alpha)(1-\theta)+\beta}}.
            \qedhere
        \]

    \end{proof}

\begin{proof} [Proof of \cref{thm::FK3K4}]
    To prove the upper bound on $f_{F,K_3}(n)$, we use the following $K_3$-free graph given by Alon in~\cite{alon1994explicit}.
    \begin{theorem}[Theorem 2.1 in~\cite{alon1994explicit}]\label{k3graph}
	If $\l$ is not divisible by $3$ and $m=2^{3\l}$, then there exists a triangle-free $d_\l$-regular graph $G_\l$ on $m$ vertices, where $d_\l=2^{\l-1}(2^{\l-1}-1)$, such that every eigenvalue $\nu$ of $G_\l$, besides the largest, satisfies
    $$
    -9\cdot2^\l-3\cdot 2^{\l/2}-1/4 \le \nu \le 4\cdot 2^\l+2\cdot 2^{\l/2}+1/4.
    $$
    \end{theorem}
    \noindent
    By the Expander Mixing Lemma (see for example~Corollary 9.2.6 in~\cite{alon2016probabilistic}), we have that every set $X \subseteq V(G_\l)$ with $|X| \ge 10^{10} \cdot 2^{2\l}$ spans at least $|X|^2/(10^{10} \cdot 2^\l)$ edges. Hence, we can apply \cref{thm::generalfFG} with $\beta = 1/3$ and $\theta = 2/3$ and get that there is a $K_3$-free graph $G'_\l$ with
    $$
    n = \Omega\left(m^{\frac{2-\alpha}{3-3\alpha}} (\log m)^{\frac{3}{1-\alpha}} \right)
    $$
    vertices and
    $$
    \alpha_{F}(G'_\l) = O\left(m^{\frac{1}{3-3\alpha}} (\log m)^{\frac{3}{1-\alpha}}\right) = O\left( n^{\frac{1}{2-\alpha}} (\log n)^{\frac{3}{2-\alpha}} \right).
    $$
    Therefore, we have
    $$
        f_{F, K_3} (n) = O\left( n^{\frac{1}{2-\alpha}} (\log n)^{\frac{3}{2-\alpha}} \right).
    $$

    For $f_{F,K_4}(n)$, we use the following construction from~\cite{mattheus2024asymptotics}.
    \begin{theorem}[Theorem 3 in~\cite{mattheus2024asymptotics}] \label{supersaturation-K4}
	For every prime power $q\ge 2^{40}$, there exists a $K_4$-free graph $G_q$ with $m=q^2(q^2-q+1)$ vertices such that for every set $X$ of at least $2^{24}q^2$ vertices of $G_q$, we have $e(G_q[X])\ge {|X|^2}/({256q})$.
    \end{theorem}
    \noindent
    We apply \cref{thm::generalfFG} with $\beta = 1/4$ and $\theta = 1/2$ and get that there is a $K_4$-free graph $G_q'$ with
    $$
        n = \Omega\left(m^{\frac{3-2\alpha}{4-4\alpha}} (\log m)^{\frac{3}{1-\alpha}} \right)
    $$
    vertices and
    $$
        \alpha_F(G_q') = O\left( m^{\frac{1}{4-4\alpha}} (\log m)^{\frac{3}{1-\alpha}}  \right)
        =O\left(n^{\frac{1}{3-2\alpha}} (\log n)^{\frac{6}{3-2\alpha}} \right).
    $$
    By Bertrand's postulate that there is always a prime number between $a$ and $2a$ for every positive integer $a$, we conclude that
    \begin{equation*}
        f_{F,K_4}(n) = O\left(n^{\frac{1}{3-2\alpha}} (\log n)^{\frac{6}{3-2\alpha}} \right). \qedhere
    \end{equation*}
\end{proof}

\begin{proof} [Proof of \cref{thm::c2kpseu}]
    By the Expander Mixing Lemma, we have that there is $\gamma = \gamma(t) > 0$ such that for every $A \subseteq V(G_m)$ with $|A| \ge \gamma m^{1/2 + 1/(4t-6)}$, we have $e(A) \ge |A|^2 \cdot \Omega(m^{-1/(2t-3)})$. Note that
    $$
        \frac{1}{2t-3} < \left(\frac{1}{2} + \frac{1}{4t-6}\right) \left(1-\alpha\right)
    $$
    by our assumption.
    The result follows by applying \cref{thm::generalfFG} with $\beta = 1/(2t-3)$, and $\theta = 1/2 + 1/(4t-6)$.
\end{proof}

\section{\texorpdfstring{Proof of \cref{thm::Ks1srVsKr+2}}{Proof of Theorem \ref{thm::Ks1srVsKr+2}}} \label{sec::boundByJanzerSudakov}

Following~\cite{janzer2023improved},
we start from the following graph, which plays a crucial role in the proof of the lower bound on $r(4,t)$~\cite{mattheus2024asymptotics} (see also~\cite{onan1972automorphisms}).
\begin{lemma} \label{lem::graphMV}
For every prime  number $q$, there is a bipartite graph $F$ with vertex bipartition $(X,Y)$ such that
\begin{enumerate}[(a)]
    \item $|X|=q^4-q^3+q^2$, $|Y|=q^3+1$,
    \item $d_F(x)=q+1$ for every $x \in X$ and $d_F(y)=q^2$ for every $y \in Y$,
    \item $F$ is $C_4$-free,
    \item there are no $\{x_k\}_{1\le k \le 4} \subseteq X$ and $\{y_{ij}\}_{1\le i < j \le 4} \subseteq Y$ such that $\{x_i, y_{ij}\}$ and $\{x_j, y_{ij}\}$ are edges in $F$ for every $1 \le i < j \le 4$, i.e.,
    $F$ does not contain the subdivision of $K_4$ as a subgraph with the four vertices of $K_4$ in $X$.
\end{enumerate}
\end{lemma}
\noindent
Let $r \ge 2$ be a fixed integer and $F$ be the graph from \cref{lem::graphMV}. For every vertex $y \in Y$, we randomly partition $N_F(y)$ into $r$ parts as $S_1(y)\cup \ldots\cup S_{r}(y)$.
Let $H$ be the (random) graph on vertex set $X$, where $\{x_1, x_2\}$ is an edge if and only if there exist a vertex $y \in Y$ and integers $1 \le i \neq j \le r$ such that $x_1 \in S_i(y)$ and $x_2 \in S_j(y)$, i.e., for every $y \in Y$, we embed a complete $r$-partite graph in $N_F(y)$. Based on (c) and (d) in \cref{lem::graphMV}, the following is proved in~\cite{janzer2023improved}. See Lemma~2.2 and the discussion after it in~\cite{janzer2023improved}.

\begin{claim} \label{cla::Kr+2free}
$H$ is always $K_{r+2}$-free.
\end{claim}
As in the original proof of the lower bound on $r(4,t)$ in~\cite{mattheus2024asymptotics}, the graph $H$ itself is not sufficient for proving \cref{thm::Ks1srVsKr+2}, as $\alpha_{K_{s_1,\ldots, s_r}}(H)$ may be too large.
We need to do some further operations on $H$.

Recall that $s = \sum_{i=1}^r s_i$.
We will consider the
$s$-uniform hypergraph $\mH$ on vertex set $X$, whose hyperedges consist of all the copies of $K_{s_1,\ldots, s_r}$ in $H$.
By definition, every $K_{s_1,\ldots, s_r}$-independent set in $H$ is an independent set in $\mH$, the number of which can be bounded by the following version of Hypergraph Container Lemma~\cite{balogh2015independent, saxton2015hypergraph}. See Corollary 2.8 in~\cite{janzer2023improved} for its proof.

Hereinafter, we denote $\{1,2,\ldots, m\}$ by $[m]$ for every positive integer $m$.
For an $s$-uniform hypergraph $\G$ and $\l \in [s]$, let $\Delta_\l (\G)$ be the \emph{maximum $\l$-degree} of $\G$, i.e., the maximum number of hyperedges containing an arbitrary fixed set of $\l$ vertices.

\begin{lemma} \label{lem::hypContainerLemma}
For every positive integer $s$ and positive real numbers $p$ and $\lambda$, the following holds. Suppose that $\G$ is an $s$-uniform hypergraph such that $p v(\G)$ and $v(\G) / \lambda$ are integers, and for every $\l \in [s]$,
\begin{equation} \label{equ::codContainerLemma}
\Delta_{\ell}(\G) \leq \lambda \cdot p^{\ell-1} \frac{e(\G)}{v(\G)}.
\end{equation}
Then, there exists a collection $\C$ of at most $v(\G)^{spv(\G)}$ sets of size at most $\left(1-\delta \lambda^{-1}\right) v(\G)$ such that for every independent set $I$ in $\G$, there exists some $R \in \C$ with $I \subseteq R$, where $\delta=2^{-s(s+1)}$.
\end{lemma}
In order to guarantee that~\eqref{equ::codContainerLemma} in \cref{lem::hypContainerLemma} is satisfied with a small value of $p$, we use the method from~\cite{janzer2023improved}, choosing a set of copies of $K_{s_1,\ldots, s_r}$ in $H$ by the following lemma.

\begin{lemma}[Lemma 2.4 in~\cite{janzer2023improved}] \label{lem::goody}
Let $q$ be a sufficiently large prime number. Then, with positive probability, for every $U \subseteq X$ such that $|U| \ge 500 r^2 q^2$, there exists some $\gamma \ge |U| / q^2$ such that the number of $y \in Y$ with $\gamma /(10 r) \le$ $\left|S_i(y) \cap U\right| \le \gamma$ for every $i \in [r]$ is at least $|U| q /(8(\log q) \gamma)$.
\end{lemma}
\noindent
We say that an instance of $H$ is \emph{good} if it satisfies the conclusion of \cref{lem::goody}.

\begin{lemma} \label{lem::numIndset}
    Let $q$ be a sufficiently large prime number and $a = q^{2- 1/(s-1)}(\log q)^3$. If $H$ is good, then the number of sets $A \subseteq X$, where $|A| = a$ and $H[A]$ is $K_{s_1,\ldots,s_r}$-free, is at most
    $$
        \left(q^{1/(s-1)}\right) ^a.
    $$
\end{lemma}
\begin{proof}
    Let $\mH$ be the $s$-uniform hypergraph on vertex set $X$, whose hyperedges consist of all the copies of $K_{s_1,\ldots, s_r}$ in $H$.
    \begin{claim} \label{cla::containerOnQSquareVertices}
        For every $U \subseteq X$ with $|U| \ge 500 s^2q^2$, there is a collection $\C_U \subseteq \mP(U)$ of at most $q^{4sq^{2- 1/(s-1)}}$ sets of size at most $(1- \Omega(\log q)^{-1})|U|$ such that for every independent set $I$ in $\mH[U]$, there is some $R \in \C_U$  with $I \subseteq R$.
    \end{claim}
    \begin{proof}
        By \cref{lem::goody}, there exist a real number $\gamma \ge |U| / q^2$ and
        a set $B \subseteq Y$ of vertices such that $|B| \ge |U|q/ (8(\log q) \gamma)$ and $\gamma /(10 r) \le$ $\left|S_i(y) \cap U\right| \le \gamma$ for every $y \in B$ and  $i \in [r]$.
        Let $\G$ be the $s$-uniform hypergraph on vertex set $U$, where $\{x_{ij}: 1\le i \le r, 1\le j \le s_i\}$ is a hyperedge in $\G$ if and only if there is $y \in B$ such that $x_{ij} \in S_i(y)$ for every $i \in [r]$ and $j \in [s_i]$.

        To apply \cref{lem::hypContainerLemma} to $\G$, we need to choose proper values for $\lambda, p$ and verify that~\eqref{equ::codContainerLemma} is satisfied.
        Note that by (c) in \cref{lem::graphMV}, we have that $|N_F(y_1) \cap N_F(y_2)| \le 1$ for distinct $y_1,y_2 \in B$. Hence, by the definition of $\G$,
        we have
        $$
            e(\G) \ge |B| \cdot \prod_{i=1}^r \b{|S_i(y) \cap U|}{s_i}
            = \frac{|U|q}{8(\log q)\gamma} \cdot \Omega\left( \left(\frac{\gamma}{10r}\right)^s \right) = \Omega \left( \frac{|U|q\gamma^{s-1}}{\log q} \right),
        $$
        which implies that
        $$
            \frac{e(\G)}{v(\G)} = \Omega \left( \frac{q\gamma^{s-1}}{\log q} \right).
        $$
        For every $x \in U$, there can be at most $d_F(x) = q+1$ vertices $y \in B$ such that $x \in N_F(y)$, so
        $$
            \Delta_1(\G) = (q+1) \cdot O\left( \gamma^{s-1} \right).
        $$
        For every set of $\l$ $(2\le \l \le s)$ vertices $\{x_1,\ldots, x_\l\} \subseteq U$, there can be at most one vertex $y \in B$ such that $\{x_1,\ldots, x_\l\} \subseteq N_F(y)$, so
        $$
            \Delta_\l(G) = O\left( \gamma^{s-\l} \right).
        $$

        Let $C$ be a sufficiently large constant compared to $s$ and $r$. Let $\lambda$ be a real number between $C \log q$ and $2C \log q$ such that $v(\G) / \lambda $ is an integer.
        Let $p$ be a real number between $(\gamma q^{1/(s-1)})^{-1}$ and $2(\gamma q^{1/(s-1)})^{-1}$ such that $pv(\G)$ is an integer. Then,
        \eqref{equ::codContainerLemma} is satisfied. By \cref{lem::hypContainerLemma}, there exists a collection $\C_U$ of at most
        $$
        v(G)^{spv(G)} \le |U|^{s \cdot 2(\gamma q^{1/(s-1)})^{-1} \cdot |U|}
        \le (q^4)^{s \cdot 2\left(\frac{|U|}{q^2} \cdot q^{1/(s-1)}\right)^{-1} \cdot |U|}
        \le
        q^{8s  q^{2- 1/(s-1)}}
        $$
        sets of size at most
        $$
            \left( 1- 2^{-s(s+1)} \lambda^{-1} \right)v(\G) = \left( 1- \Omega\left((\log q)^{-1}\right)\right) |U|
        $$
        such that for every independent set $I$ in $\G$, there exists some $R \in \C_U$  with $I \subseteq R$. The conclusion in \cref{cla::containerOnQSquareVertices} follows by noting that $\G$ is a subhypergraph of $\mH$, so every independent set in $\mH[U]$ is an independent set in $\G$.
    \end{proof}

    \begin{claim} \label{cla::containerOnTheWholeGraph}
    There is a collection $\C \subseteq \mP(X) $ of at most $q^{O(q^{2-1 /(s-1)}(\log q)^2)}$ sets of size at most $500 s^2 q^2$ such that for every independent set $I$ in $\mH$, there is some $R \in \C$ such that $I \subseteq R$.
    \end{claim}
    \begin{proof}
        We prove by induction that for every positive integer $j$, there is a collection $\C^{(j)}$ of at most  $q^{8jsq^{2-1/(s-1)}}$ sets of size at most $\max\left\{500s^2q^2,\,(1-\Omega((\log q)^{-1}))^j|X|\right\}$ such that for every independent set $I$ in $\mH$, there is some $R \in \C^{(j)}$ such that $I \subseteq R$. \cref{cla::containerOnTheWholeGraph} then follows by letting $j$ be a suitable integer of order $\Theta((\log q)^2)$.

        The case $j=1$ holds by applying \cref{cla::containerOnQSquareVertices} with $U = X$. Now, suppose we have a proper collection $\C^{(j)}$. For every $U \in \C^{(j)}$ with $|U| > 500s^2q^2$, we apply \cref{cla::containerOnQSquareVertices} to $U$ and get a collection $\C_U$. Let
        $$
            \C^{(j+1)} \ce \left\{U \in \C^{(j)} : |U| \le 500s^2q^2 \right\} \cup \bigcup_{U \in \C^{(j)}: |U| > 500s^2q^2} \C_U.
        $$
        By the induction hypothesis, we have
        $$
            \C^{(j+1)} \le |\C^{(j)}| \cdot q^{8sq^{2-1/(s-1)}} \le q^{8jsq^{2-1/(s-1)}} \cdot q^{8sq^{2-1/(s-1)}} =  q^{8(j+1)sq^{2-1/(s-1)}}.
        $$
        For every set $R \in \C^{(j+1)}$, either $|R| \le 500s^2q^2$ or there is $U \in \C^{(j)}$ such that $R \in \C_U$, in which case we have
        $$
        |R| \le (1-\Omega((\log q)^{-1}))|U| \le (1-\Omega((\log q)^{-1}))^{j+1}|X|.
        $$
        So, $|R| \le \max\left\{500s^2q^2,\,(1-\Omega((\log q)^{-1}))^{j+1}|X|\right\}$. For every independent set $I \in \mH$, there is some $U \in \C^{(j)}$ such that $I \subseteq U$. If $|U| \le 500s^2q^2$, then $I\subseteq U \in \C^{(j+1)}$; if $|U| > 500s^2q^2$, then there is some $R \in \C_U$ such that $I \subseteq R \in \C_U \subseteq \C^{(j+1)}$. Therefore, $C^{(j+1)}$ indeed has the properties claimed.
    \end{proof}

    \noindent
    By \cref{cla::containerOnTheWholeGraph}, the number of independent sets of size $a = q^{2- 1/(s-1)}(\log q)^3$ in $\mH$ is at most
    \begin{equation*}
        q^{O(q^{2-1 /(s-1)}(\log q)^2)} \cdot \b{500 s^2 q^2}{a}
        \le
        q^{O(q^{2-1 /(s-1)}(\log q)^2)} \cdot \left(\frac{500e\cdot  s^2 q^2}{a}\right)^a
        \le
        \left(q^{1/(s-1)} \right)^a. \qedhere
    \end{equation*}

\end{proof}

\begin{proof}[Proof of \cref{thm::Ks1srVsKr+2}]
    Let $q$ be a sufficiently large prime number and fix a good $H$, whose existence is guaranteed by \cref{lem::goody}. Note that $H$ is $K_{r+2}$-free by \cref{cla::Kr+2free}.
    Let $a \ce q^{2- 1/(s-1)}(\log q)^3$ and $\A$ be the collection of sets $A \subseteq X$, where $|A| = a$ and $H[A]$ is $K_{s_1,\ldots,s_r}$-free. By \cref{lem::numIndset}, we have $|\A| \le \left(q^{1/(s-1)}\right) ^a$.

    We choose a random subset $X'$ of $X$, where every vertex in $X$ is kept independently with probability $q^{-1/(s-1)}$. Let $H' \ce H[X']$ and $\A'$ be the (random) collection of sets $\{A \in \A: A\subseteq X'\}$. Let $H''$ be the (random) graph obtained from $H'$ by removing one vertex in every $A \in \A'$. Note that $H''$ is also $K_{r+2}$-free and does not contain any $K_{s_1,\ldots, s_r}$-free set of size $a$. The expected number of vertices in $H''$ is at least
    $$
        \E[|X'| - |\A'|] = \E[|X'|] - \E[|\A'|] \ge |X| \cdot q^{-1/(s-1)} - 1 \ge \frac{1}{2}q^{4- 1/(s-1)}.
    $$
    Therefore, for every sufficiently large prime number $q$, there exists a $K_{r+2}$-free graph $H''$ with at least $\frac{1}{2}q^{4- 1/(s-1)}$ vertices and $\alpha_{K_{s_1,\ldots, s_r}}(G) < a = q^{2- 1/(s-1)}(\log q)^3$.
    By Bertrand's postulate that there is always a prime number between $a$ and $2a$ for every positive integer $a$, we conclude that
    \begin{equation*}
        f_{K_{s_1,\ldots,s_r}, K_{r+2}}(n) = O\left( n^{\frac{2s -3}{4s -5}} (\log n)^{3} \right). \qedhere
    \end{equation*}
\end{proof}

\section{\texorpdfstring{Remark on $\ex(G(n,p), C_4)$}{Remark on ex(G(n,p),C4)}} \label{sec::GnpC4}

The previous bounds on $\ex(G(n,p), C_4)$ are the   following.
\begin{theorem}[\cite{furedi1994random, kohayakawa1998extremal, morris2016number}]
For every $\delta>0$, we have
\begin{alignat*}{3}
    & \ex(G(n,p), C_4) = \Theta\left( p^{1/2} n^{3/2} \right) \quad && \textrm{ if } \quad p\ge n^{-1/3} (\log n)^{4},&& \\
    & O\left(n^{4/3} (\log n)^2 \right) = \ex(G(n,p),C_4) = \Omega\left(n^{4/3}(\log n)^{1/3}\right) \quad && \textrm{ if } \quad n^{-1/3}(\log n)^4 > p> n^{-1/3-\delta},&& \\
    & \ex(G(n,p),C_4) = \Theta\left(n^{4/3} (\log (pn^{2/3}))^{1/3}\right) \quad && \textrm{ if } \quad n^{-1/3-\delta} \ge p > 2n^{-2/3},&& \\
    & \ex(G(n,p),C_4) = \Theta\left( p\b{n}{2} \right) \quad && \textrm{ if } \quad 2n^{-2/3} \ge  p> 10n^{-2}.&&
\end{alignat*}
\end{theorem}
\noindent
\cref{pro::impBoundGnpC4} improves the bounds in the first two ranges as follows:
\begin{alignat*}{3}
    & \ex(G(n,p), C_4) = \Theta\left( p^{1/2} n^{3/2} \right) \quad && \textrm{ if } \quad p \ge n^{-1/3} (\log n)^{8/3},&& \\
    & O\left(n^{4/3} (\log n)^{4/3} \right) = \ex(G(n,p),C_4) = \Omega\left(n^{4/3}(\log n)^{1/3}\right) \quad && \textrm{ if } \quad n^{-1/3}(\log n)^{8/3} > p> n^{-1/3-\delta}.
\end{alignat*}
Note that the ``construction'' for the lower bound is the same, and we extend the range of $p$, where it is best possible. In order to obtain the claimed improvement, we use the following theorem and the first moment method.
Denote $f(C_4, n, T)$ the number of (labelled) $C_4$-free graphs on $n$ vertices with $T$ edges.

\begin{theorem}[Lemma 5 in~\cite{kohayakawa1998extremal}]\label{KKS}
	Let $c=(2^{65}/3)^{1/3}$ and $\eps=3/(13\cdot 2^{65})$. For every integer $T > c n^{4/3} (\log n)^{1/3}$, we have
    \[
    f(C_4, n, T)<\left(\frac{4en^3}{T^2}\right)^T\exp \left(\eps^{-1}\left(n\log n\right)^{4/3}\right).
    \]
\end{theorem}

\begin{proof}[Proof of \cref{pro::impBoundGnpC4}]
    Let $C$ be a sufficiently large real number and $T$ be an integer between $Cp^{1/2}n^{3/2}$ and $2Cp^{1/2}n^{3/2}$. Let $c=(2^{65}/3)^{1/3}$ and $\eps=3/(13\cdot 2^{65})$, as in \cref{KKS}.
    Note that
    \begin{equation} \label{equ::lowBoundT}
    T \ge Cp^{1/2} n^{3/2} \ge C n^{-1/6}(\log n)^{4/3} \cdot n^{3/2}  = Cn^{4/3} (\log n)^{4/3} > c n^{4/3}(\log n)^{1/3}.
    \end{equation}
    So, by \cref{KKS}, the expected number of $C_4$-free subgraphs with $T$ edges in $G(n,p)$ is at most
    \begin{align*}
        p^T\left(\frac{4en^3}{T^2}\right)^T\exp \left(\eps^{-1}\left(n\log n\right)^{4/3}\right)
        &\le \left(\frac{4epn^3}{C^2 p n^3}\right)^T\exp \left(\eps^{-1}\left(n\log n\right)^{4/3}\right) \\
        &= \left(\frac{4e}{C^2}\right)^T\exp \left(\eps^{-1}\left(n\log n\right)^{4/3}\right)
        = o(1).
    \end{align*}
    Therefore, with high probability, there is no $C_4$-free subgraph with $T$ edges in $G(n,p)$ and hence, we have $\ex(G(n,p), C_4) <  T$.
\end{proof}

\section*{Acknowledgment}
We thank the anonymous referees for many helpful comments that improved the quality of the manuscript.
The authors are grateful to Noga Alon,
Sam Mattheus, Dhruv Mubayi, Jacques Verstra{\"e}te, and Ethan White for helpful discussions.

\end{document}